\documentclass[11pt]{article}
\usepackage{amsmath}
\usepackage{amsfonts}
\usepackage{amssymb}
\usepackage{hyperref}
\usepackage{indentfirst}
\usepackage{mathrsfs}

\usepackage[top=1in,bottom=1in,left=1.25in,right=1.25in]{geometry}
\date{}
\allowdisplaybreaks

\begin{document}

\renewcommand{\baselinestretch}{1.2}
\renewcommand{\arraystretch}{1.0}

\title{\bf Another construction of the braided $T$-category}
\author
{
  \textbf{Tao Yang} \footnote{College of Science, Nanjing Agricultural University, Nanjing 210095, Jiangsu, CHINA.
             E-mail: tao.yang.seu@gmail.com}
   , \quad
  \textbf{Xiaoyan Zhou}
  \footnote{Corresponding author. College of Science, Nanjing Agricultural University, Nanjing 210095, Jiangsu, CHINA.
            E-mail: zhouxy@njau.edu.cn}
}

\maketitle

\begin{center}
\begin{minipage}{12.cm}

 {\bf Abstract } This paper introduces group-cograded monoidal Hom-Hopf algebras, and shows that this kind of
 group-cograded monoidal Hom-Hopf algebras are monoidal Hom-Hopf algebras in the Turaev category $\mathcal{J}_{k}$
 introduced by Canepeel and De Lombaerde.
 Then we define the $p$-Yetter-Drinfeld category over a group-cograded monoidal Hom-Hopf algebra,
 and construct a new kind of braided $T$-categories.
\\

 {\bf Key words } Monoidal Hom-Hopf algebra, group-cograded, Yetter-Drinfeld module,  Turaev category.
\\

 {\bf Mathematics Subject Classification:} 16T05

\end{minipage}
\end{center}
\normalsize

\section{Introduction}
\def\theequation{\thesection.\arabic{equation}}
\setcounter{equation}{0}

 Turaev (see \cite{T00,T08}) generalized quantum invariants of $3$-manifolds to the case of a $3$-manifold $M$
 endowed with a homotopy class of maps $M\longrightarrow  K(G, 1)$,
 where $K(G, 1)$ is the Eilenberg-MacLane space as a target space determined by a group $G$.
 Braided crossed categories based on a group $G$ (i.e. braided $T$-categories in the sense of Turaev,
 is braided monoidal categories in Freyd-Yetter categories of crossed $G$-sets (see \cite{FY89})
 play a key role in the construction of these homotopy invariants.

 One main question naturally arises: how to construct classes of new braided $T$-categories?
 For this, Panaite and Staic (see \cite{PS07}) and Zunino (see \cite{Z04}) gave some interesting constructions in Hopf algebra case.
 Following their motivated idea, Yang et al (see \cite{YW11,YZM13,ZY13}) generalized their work to
 multiplier Hopf algebra and weak Hopf algebra cases, also the authors get some new classes of braided $T$-categories.

 In \cite{CG11}, Caenepeel and Goyvaerts illustrated Hom-structures from the point of view of monoidal categories and introduced monoidal
 Hom-Hopf algebras. Many of the classical results in Hopf algebra theory can be generalized to the monoidal Hom-Hopf algebra setting.
 Such as in \cite{YW14}, You and Wang generalized the braided $T$-category construction of Panaite and Staic to the monoidal Hom-Hopf algebra case.
 Motivated by this, in this paper we mainly generalize the constructions shown in \cite{Z04}, replacing their Hopf group-coalgebras
 by so-called monoidal Hom-Hopf group-coalgebras, and provide new examples of braided $T$-categories.

 The paper is organized in the following way.

 In section 2, we recall some notions which we will use in the following, such as Hom-category, monoidal Hom-Hopf algebras,
 and their modules and comodules.

 In section 3, we give the definitions of group-cograded monoidal Hom-Hopf algebras, monoidal Hom-Hopf $T$-coalgebras,
 and show some examples.

 In section 4, based on the $p$-Yetter-Drinfeld module category over a monoidal Hom-Hopf $T$-coalgebra,
 we give a method to construct a new class of braided $T$-categories, which is the main result of this paper.

\section{Preliminaries}
\def\theequation{\thesection.\arabic{equation}}
\setcounter{equation}{0}

 Throughout this paper, all vector spaces, tensor products, and morphisms are over a fixed field $k$.
 We use the Sweedler's notation for terminologies on coalgebras. For a coalgebra $C$, we write comultiplication
 $\Delta(c)=c_{(1)} \otimes c_{(2)}$ for any $c \in C$, in which we often omit the summation symbols for convenience.
 And we denote $Id_{M}$ for the identity map from $M$ to $M$.

 \subsection{Monoidal Hom-Hopf algebras}

 Let $\mathcal{M}_{k} = (\mathcal{M}_{k}, \otimes, k, a, l, r)$ be the category of $k$-modules.
 There is a new monoidal category $\mathcal{H}(\mathcal{M}_{k})$:
 the objects are pairs $(M, \mu)$, where $M\in \mathcal{M}_{k}$ and $\mu \in Aut_{k}(M)$,
 and the morphisms $f: (M, \mu)\longrightarrow (N, \nu)$ in $\mathcal{M}_{k}$ such that $\nu\circ f=f \circ \mu$.
 For any objects  $(M, \mu)$ and $(N, \nu)$, the monoidal structure is given by
 \begin{eqnarray*}
 (M, \mu) \otimes (N, \nu) = (M\otimes N, \mu\otimes \nu) \qquad \mbox{and} \qquad (k, Id_{k}).
 \end{eqnarray*}

 Briefly speaking, all Hom-structures are objects in the monoidal category $\mathcal{\widetilde{H}}(\mathcal{M}_{k})$
 = $(\mathcal{H}(\mathcal{M}_{k}), \otimes, (k, Id), \widetilde{a}, \widetilde{l}, \widetilde{r})$, where the associator $\widetilde{a}$
 and the unitors $\widetilde{l}, \widetilde{r}$ are given by
 \begin{eqnarray*}
 \widetilde{a}_{M, N, L} = {a}_{M, N, L} \circ ((\mu \otimes Id_{N}) \otimes \varsigma^{-1})
                         = (\mu \otimes (Id_{N} \otimes \varsigma^{-1}))\circ {a}_{M, N, L},\\
 \widetilde{l}_{M} = \mu \circ l_{M} = l_{M} \circ (Id \otimes \mu), \qquad
 \widetilde{r}_{M} = \mu \circ r_{M} = r_{M} \circ (\mu \otimes Id),
 \end{eqnarray*}
 for any $(M, \mu), (N, \nu), (L, \varsigma) \in \mathcal{\widetilde{H}}(\mathcal{M}_{k})$ (see \cite{CG11}).
 The category $\mathcal{\widetilde{H}}(\mathcal{M}_{k})$ is called the Hom-category associated to the monoidal category $\mathcal{M}_{k}$.
 \\

 In the following, we will recall from \cite{CG11} some definitions about Hom-structures.

 A \emph{unital monoidal Hom-associative algebra} $(A, \alpha)$ is an object in the Hom-category $\mathcal{\widetilde{H}}(\mathcal{M}_{k})$,
 together with an element $1_{A}\in A$ and a linear map $m: A\otimes A\longrightarrow A, a\otimes b\mapsto ab$ such that
 for all $a, b, c \in A$,
 \begin{eqnarray}
 \alpha(a)(bc)=(ab)\alpha(c), \quad a1_{A} = 1_{A}a = \alpha(a), \\
 \alpha(ab)=\alpha(a)\alpha(b), \qquad \alpha(1_{A})=1_{A}.
 \end{eqnarray}
 In this paper, the algebras we mainly considered are this kind of unital monoidal Hom-associative algebras,
 and in the follwing we call them the monoidal Hom-algebras without confusion.
 Let $(A, \alpha)$ and $(A', \alpha')$ be two monoidal Hom-algebras.
 A Hom-algebra map $f: (A, \alpha)\longrightarrow(A', \alpha')$ is a linear map such that $f\circ\alpha=\alpha'\circ f$,
 $f(ab)=f(a)f(b)$ and $f(1_{A})=1_{A'}$.

 A \emph{counital monoidal Hom-coassociative coalgebra} (monoidal Hom-coalgebra in short) $(C, \gamma)$
 is an object in $\mathcal{\widetilde{H}}(\mathcal{M}_{k})$ together with linear maps
 $\Delta: C\longrightarrow C\otimes C, c\mapsto c_{(1)}\otimes c_{(2)}$
 and $\varepsilon: C\longrightarrow k$ such that for any $c\in C$,
 \begin{eqnarray}
 \gamma^{-1}(c_{(1)}) \otimes \Delta(c_{(2)}) = \Delta(c_{(1)}) \otimes \gamma^{-1}(c_{(2)}),
 c_{(1)}\varepsilon(c_{(2)})=\varepsilon(c_{(1)})c_{(2)}=\gamma^{-1}(c), \\
 \Delta(\gamma^{-1}(c))=\gamma^{-1}(c_{(1)}) \otimes \gamma^{-1}(c_{(2)}), \qquad
 \varepsilon(\gamma^{-1}(c))=\varepsilon(c).
 \end{eqnarray}
 Let $(C, \gamma)$ and $(C', \gamma')$ be two monoidal Hom-coalgebras. A Hom-coalgebra map $f: (C, \gamma)\longrightarrow(C', \gamma')$
 is a linear map such that $f\circ \gamma = \gamma'\circ f$, $\Delta'\circ f = (f\otimes f)\Delta$ and $\varepsilon'\circ f= \varepsilon$.

 Recall from \cite{LS14} that a \emph{monoidal Hom-bialgebra} $H=(H, \alpha, m, 1_{H}, \Delta, \varepsilon)$
 is a bialgebra in $\mathcal{\widetilde{H}}(\mathcal{M}_{k})$.
 This means $(H, \alpha, m, 1_{H})$ is a monoidal Hom-algebra and $(H, \alpha, \Delta, \varepsilon)$ is a monoidal Hom-coalgebra such that
 $\Delta, \varepsilon$ are Hom-algebra maps, i.e., for any $h, g\in H$,
 \begin{eqnarray}
 \Delta(hg) = \Delta(h)\Delta(g), \Delta(1_{H}) = 1_{H} \otimes 1_{H}, \\
 \varepsilon(hg)=\varepsilon(h)\varepsilon(g), \qquad \varepsilon(1_{H})=1_{k}.
 \end{eqnarray}
 A monoidal Hom-bialgebra $H=(H, \alpha)$ is called a \emph{monoidal Hom-Hopf algebra} if there is a morphism (called the antipode)
 $S: H\longrightarrow H$ in $\mathcal{\widetilde{H}}(\mathcal{M}_{k})$ (i.e., $S\circ \alpha = \alpha\circ S$),
 which is the convolution inverse of the identity morphism $Id_{H}$, this means for any $h\in H$
 \begin{eqnarray}
 S(h_{(1)})h_{(2)} = \varepsilon(h)1_{H} = h_{(1)}S(h_{(2)}).
 \end{eqnarray}

 \subsection{Actions and coactions over monoidal Hom-(co)algebras}

 Let $(A, \alpha)$ be a monoidal Hom-algebra.
 A \emph{left $(A, \alpha)$-Hom-module} consists of $(M, \mu)\in \mathcal{\widetilde{H}}(\mathcal{M}_{k})$
 with a morphism $\psi: A\otimes M\longrightarrow M, a\otimes m \mapsto am$ such that for all $a, b\in A$ and $m\in M$,
 \begin{eqnarray}
 \alpha(a)(bm)=(ab)\mu(m), \quad 1_{A}m= \mu(m),\\
 \mu(am)=\alpha(a)\mu(m).
 \end{eqnarray}
 Let $(M, \mu)$ and $(N, \nu)$ be two left $(A, \alpha)$-Hom-modules,
 a morphism $f: M\longrightarrow N$ is called left $A$-linear if for any $a\in A, m\in M$, $f(am)=af(m)$ and $f\circ\mu = \nu\circ f$.
 \\

 Let $(C, \gamma)$ be a monoidal Hom-coalgebra.
 A right $(C, \gamma)$-Hom-comodule  is an object $(M, \mu) \in \mathcal{\widetilde{H}}(\mathcal{M}_{k})$ together with a linear map
 $\rho: M\longrightarrow M\otimes C, m\mapsto m_{(0)}\otimes m_{(1)}$ in $\mathcal{\widetilde{H}}(\mathcal{M}_{k})$ such that
 \begin{eqnarray}
 \mu^{-1}(m_{(0)}) \otimes \Delta(m_{(1)}) = m_{(0)(0)} \otimes (m_{(0)(1)} \otimes \gamma^{-1}(m_{(1)})),
 \quad m_{(0)}\varepsilon(m_{(1)}) = \mu^{-1}(m),
 \end{eqnarray}
 for all $m\in M$. Indeed, $\rho\in \mathcal{\widetilde{H}}(\mathcal{M}_{k})$ means
 \begin{eqnarray}
 \rho(\mu(m)) = \mu(m_{(0)}) \otimes \gamma(m_{(1)}).
 \end{eqnarray}
 Let $(M, \mu)$ and $(N, \nu)$ be two right $(C, \gamma)$-Hom-comodules,
 a morphism $g: M\longrightarrow N$ is called right $C$-colinear if $g\circ\mu = \nu\circ g$ and for any $m\in M$,
 $g(m)_{(0)} \otimes g(m)_{(1)} = g(m_{(0)}) \otimes m_{(1)}$.

\section{Group-cograded monoidal Hom-Hopf algebras}
\def\theequation{\thesection.\arabic{equation}}
\setcounter{equation}{0}

 By Proposition 2.5 in paper \cite{CK06}, Hopf group-coalgebras (or in the other words group-cograded Hopf algebras)
 are Hopf algebras in Turaev category $\mathcal{J}_{k}$.
 It turns out that many of the classical results in Hopf algebra theory can be generalized to the Hopf group-coalgebra setting.
 Motivated by this, in the following, we will introduce the group-cograded monoidal Hom-Hopf algebras.
 Firstly, we show the definition of group-cograded Hom-coalgebra. Let $G$ be a group with unit $e$.
 \\

 \textbf{Definition \thesection.1}
 A group-cograded monoidal Hom-coalgebra $(C, \gamma)$ over group $G$ is a family of objects $\{(C_{p}, \gamma_{p})\}_{p\in G}$
 in $\mathcal{\widetilde{H}}(\mathcal{M}_{k})$ together with linear maps $\Delta=\{\Delta_{p, q}\}_{p, q\in G}$,
 $\Delta_{p, q}: C_{pq}\longrightarrow C_{p}\otimes C_{q}, c \mapsto c_{(1, p)}\otimes c_{(2, q)}$
 and $\varepsilon: C_{e}\longrightarrow k$ such that
 \begin{eqnarray}
 \gamma^{-1}_{p}(c_{(1, p)}) \otimes \Delta_{q, r}(c_{(2, qr)}) = \Delta_{p, q}(c_{(1, pq)}) \otimes \gamma^{-1}_{r}(c_{(2, r)}),
 && \forall c\in C_{pqr}, \\
  c_{(1, p)}\varepsilon(c_{(2, e)})=\varepsilon(c_{(1, e)})c_{(2, p)}=\gamma^{-1}_{p}(c), && \forall c\in C_{p}, \\
 \Delta_{p, q}(\gamma^{-1}_{pq}(c_{pq}))=\gamma^{-1}_{p}(c_{(1, p)}) \otimes \gamma^{-1}_{q}(c_{(2, q)}), && \forall c\in C_{pq}, \\
 \varepsilon(\gamma^{-1}_{e}(c))=\varepsilon(c), && \forall c\in C_{e}.
 \end{eqnarray}

 Let $(C, \gamma)$ and $(C', \gamma')$ be two group-cograded monoidal Hom-coalgebras over $G$.
 A Hom-coalgebra map $f: (C, \gamma)\longrightarrow(C', \gamma')$
 is a family of linear maps $\{f_{p}\}_{p\in G}$, $f_{p}: (C_{p}, \gamma_{p})\longrightarrow(C'_{p}, \gamma'_{p})$
 such that $f_{p}\circ \gamma_{p} = \gamma'_{p}\circ f_{p}$, $\Delta'_{p, q}\circ f_{pq} = (f_{p}\otimes f_{q})\Delta_{p, q}$
 and $\varepsilon' \circ f_{e}= \varepsilon$.
 \\

 \textbf{Remark \thesection.2}
 (1) There are two trival examples: (a) If group $G$ has only one element $e$,
 then this definition is exactly the usual monoidal Hom-coalgebra.
 (b) If $(D, \varrho)$ is a monoidal Hom-coalgebra, let $C_{p}=D$ and $\gamma_{p}=\varrho$ for every $p\in G$,
 then $(C, \gamma)$ is a group-cograded monoidal Hom-coalgebra.

 (2) Group-cograded monoidal Hom-coalgebras are monoidal Hom-coalgebras in Turaev category $\mathcal{J}_{k}$.
 Indeed, let $(C, \gamma)$ be a group-cograded monoidal Hom-coalgebra, (similar to the proof of Proposition 2.2 in \cite{CK06}) we define
 $(\underline{C}, \underline{\gamma})=\big(\underline{C}=(G, C), \underline{\Delta}=(m, \Delta),
 \underline{\varepsilon}=(i, \varepsilon), \underline{\gamma}=(Id_{G}, \gamma) \big)$ with the monoidal Hom-coalgebra structure
 \begin{align*}
   \underline{C}&\stackrel{\underline{\varepsilon}}{\rightarrow} \underline{k}
   & \underline{C}\stackrel{\underline{\Delta}}{\rightarrow} \underline{C}\otimes \underline{C}
   && \underline{C}\stackrel{\underline{\gamma}}{\rightarrow} \underline{C} \\
   G &\stackrel{i}{\leftarrow} \{\ast \}
   & G\stackrel{m}{\leftarrow} G\times G
   && G\stackrel{Id_{G}}{\leftarrow} G\\
   C_{e}=C_{i(\ast)}&\stackrel{\varepsilon}{\rightarrow} k
   & C_{gh}=C_{m(g, h)}\stackrel{\Delta_{g, h}}{\rightarrow} C_{g}\otimes C_{h}
   && C_{g}=C_{Id_{G}(g)}\stackrel{\gamma_{g}}{\rightarrow} C_{g}.
 \end{align*}
 Then $(\underline{C}, \underline{\gamma})$ is a monoidal Hom-coalgebra in $\mathcal{J}_{k}$.
 \\

 \textbf{Definition \thesection.3}
 A group-cograded monoidal Hom-Hopf algebra $(H=\bigoplus_{p\in G}H_{p}, \alpha = \{\alpha_{p}\}_{p\in G})$ over $G$
 is a group-cograded monoidal Hom-coalgebra
 where each $(H_{p}, \alpha_{p})$ is a monoidal Hom-algebra with multiplication $m_{p}$ and unit $1_{p}$
 endowed with antipode $S=\{S_{p}\}_{p\in G}$,
 $S_{p}: H_{p}\longrightarrow H_{p^{-1}} \in \mathcal{\widetilde{H}}(\mathcal{M}_{k})$
 (i.e., $S_{p}\circ\alpha_{p} = \alpha_{p^{-1}}\circ S_{p}$) such that
 \begin{eqnarray}
 \Delta_{p, q}(hg) = \Delta_{p, q}(h)\Delta_{p, q}(g), \quad \Delta_{p, q}(1_{pq}) = 1_{p} \otimes 1_{q}, && \forall h, g\in H_{pq}\\
 \varepsilon(hg)=\varepsilon(h)\varepsilon(g), \qquad \qquad \quad \varepsilon(1_{e})=1_{k}, && \forall h, g\in H_{e}\\
 S_{p^{-1}}(h_{(1, p^{-1})})h_{(2, p)} = \varepsilon(h)1_{p} = h_{(1, p)}S_{p^{-1}}(h_{(2, p^{-1})}), && \forall h \in H_{e}.
 \end{eqnarray}

 Note also that the $(H_{e}, \alpha_{e}, m_{e}, 1_{e}, \Delta_{e, e}, \varepsilon, S_{e})$ is a monoidal Hom-Hopf algebra
 in the usual sense of the word. We call it the neutral component of $H$.
 \\

 In the following, we give a family of examples of group-cograded monoidal Hom-Hopf algebras.

 \textbf{Example \thesection.4}
 Let $H= \big(\bigoplus_{p\in G}H_{p}, m=\{m_{p} \}_{p\in G}, 1,
 \Delta=\{\Delta_{p, q}\}_{p, q\in G}, \varepsilon, S=\{S_{p}\}_{p\in G} \big)$
 be a Hopf group-coalgebra, then for $\alpha \in Aut(H)$ with $\alpha_{p}=\alpha|_{H_{p}}: H_{p} \longrightarrow H_{p}$,
 there is a group-cograded monoidal Hom-Hopf algebra $(H, \alpha)$ with
 multiplication $\{\alpha_{p}\circ m_{p} \}_{p\in G}$ and comultiplication $\{\Delta_{p, q}\circ \alpha^{-1}_{pq}\}_{p, q\in G}$.

 Indeed, for Hopf group-coalgebra $H$, the $p$-component $(H_{p}, m_{p})$ is an algebra, by Proposition 2.3 in \cite{CG11}
 $(H_{p}, \alpha_{p}\circ m_{p})$ is a monoidal Hom-algebra. And we can easily check that $(H, \alpha)$ satisfies (3.1)-(3.7).
 Take (3.1) and (3.7) for example, denote $\overline{m}_{p}=\alpha_{p}\circ m_{p}$ and
 $\overline{\Delta}_{p, q} = \Delta_{p, q}\circ \alpha^{-1}_{pq}$, then for (3.1),
 \begin{eqnarray*}
 && (\alpha^{-1}_{p} \otimes \overline{\Delta}_{q, r}) \overline{\Delta}_{p, qr} (c) \\
 &=& (\alpha^{-1}_{p} \otimes \overline{\Delta}_{q, r}) (\alpha^{-1}_{p}(c_{(1, p)}) \otimes \alpha^{-1}_{qr}(c_{(2, qr)}) )
  = \alpha^{-2}_{p}(c_{(1, p)}) \otimes(\alpha^{-2}_{q}(c_{(21, q)}) \otimes \alpha^{-2}_{r}(c_{(22, r)})) \\
 &=& (\alpha^{-2}_{p}(c_{(11, p)}) \otimes \alpha^{-2}_{q}(c_{(12, q)}) ) \otimes \alpha^{-2}_{r}(c_{(2, r)})
  = \overline{\Delta}_{p, q}(\alpha^{-1}_{pq}(c_{(1, pq)})) \otimes \alpha^{-2}_{r}(c_{(2, r)}) \\
 &=& (\overline{\Delta}_{p, q}\otimes \alpha^{-1}_{r}) \overline{\Delta}_{pq, r}(c),
 \end{eqnarray*}
 where $c\in H_{pqr}$. For (3.7), we only check the left part, the right part is similar. For $h\in H_{e}$,
 \begin{eqnarray*}
 \overline{m}_{p} (S_{p^{-1}} \otimes Id_{H_{p}}) \overline{\Delta}_{p^{-1}, p}(h)
 &=& \overline{m}_{p} \big(S_{p^{-1}}\alpha^{-1}_{p^{-1}}(h_{(1, p^{-1})})  \otimes \alpha^{-1}_{p}(h_{(2, p)}) \big) \\
 &=& S_{p^{-1}}(h_{(1, p^{-1})}) h_{(2, p)} = \varepsilon(h)1_{p}.
 \end{eqnarray*}

 \textbf{Proposition \thesection.5}
 Group-cograded monoidal Hom-Hopf algebras are monoidal Hom-Hopf algebras in Turaev category $\mathcal{J}_{k}$,
 and also group-cograded Hopf algebra in Hom-category $\mathcal{\widetilde{H}}(\mathcal{M}_{k})$.

 \emph{Proof} Similar to the proof of Proposition 2.5 in \cite{CK06}, let $s: G \longrightarrow G$, $s(g)=g^{-1}$,
 and consider a map $\underline{S} = (s, S) : \underline{H} \longrightarrow \underline{H}$ in $\mathcal{J}_{k}$,
 we can easily get the first result.
 \\

 Let $(H, \alpha)$ be a group-cograded monoidal Hom-Hopf algebra  over $G$. If there is a group homomorphism $\pi: G\longrightarrow Aut(H)$,
 where $Aut(H)$ denotes the group of Hom-algebra automorphism on $(H, \alpha)$,
 we call $\pi$ an admissible action of $G$ on $(H, \alpha)$, if also the following requirements hold:
 \begin{enumerate}
   \item[(1)] $\Delta\circ\pi_{p} = (\pi_{p} \otimes \pi_{p})\Delta$.
   \item[(2)] $\pi_{p}(H_{q}) = H_{\vartheta_{p}(q)}, \pi_{p} \circ \alpha_{q} = \alpha_{\vartheta_{p}(q)} \circ \pi_{p}$,
              where $\vartheta$ is an action of the group $G$ on itself.
   \item[(3)] $\pi_{\vartheta_{p}(q)} = \pi_{pqp^{-1}}$, $\alpha_{\vartheta_{p}(q)} = \alpha_{pqp^{-1}}$.

              If $\vartheta$ is the adjoint action itself, $\pi$ is called a crossing.
 \end{enumerate}

 \textbf{Definition \thesection.6}
 A group-cograded monoidal Hom-Hopf algebra $(H=\bigoplus_{p\in G}H_{p}, \alpha)$ over $G$ is said to be a  monoidal Hom-Hopf $T$-coalgebra,
 provided it is endowed with a crossing $\pi$ preserves the comultiplication and counit, i.e.,
  \begin{eqnarray}
 (\pi_{q} \otimes \pi_{q})\circ\Delta_{p, r} = \Delta_{qpq^{-1}, qrq^{-1}}\circ\pi_{q}, \qquad \varepsilon\circ\pi_{q}=\varepsilon
 \end{eqnarray}
 and $\pi$ is multiplicative in the sense that $\pi_{pq}=\pi_{p}\pi_{q}$ for all $p, q, r \in G$.
 \\

 \textbf{Proposition \thesection.7}
 Monoidal Hom-Hopf $T$-coalgebras are Hopf $T$-coalgebras in Hom-category $\mathcal{\widetilde{H}}(\mathcal{M}_{k})$.
 \\

 \textbf{Example \thesection.8}
 Let $H=(\bigoplus_{p\in G}H_{p}, \{m_{p} \}_{p\in G}, 1, \{\Delta_{p, q}\}_{p, q\in G}, \varepsilon, \{S_{p}\}_{p\in G})$
 be a Hopf $T$-coalgebra with the crossing action $\pi$.
 Then there is a group-cograded monoidal Hom $T$-coalgebra $(H, \pi)$ with
 multiplication $\{\pi_{p}\circ m_{p} \}_{p\in G}$ and comultiplication $\{\Delta_{p, q}\circ \pi^{-1}_{pq}\}_{p, q\in G}$.

\section{A new class of braided $T$-categories}
\def\theequation{\thesection.\arabic{equation}}
\setcounter{equation} {0}

 In this section, we introduce the definition of $p$-Yetter-Drinfeld modules over a monoidal Hom-Hopf $T$-coalgebra $(H, \alpha)$
 with bijective antipode,
 and show a method to construct the category $\mathcal{YD}(H)$ of Yetter-Drinfeld modules over $(H, \alpha)$, which is a braided $T$-category.
 \\

 \textbf{Definition \thesection.1}
 Let $(H=\bigoplus_{p\in G}H_{p}, \alpha = \{\alpha_{p}\}_{p\in G})$ be a monoidal Hom-Hopf $T$-coalgebra over group $G$
 and $p$ a fixed element in $G$.
 A (left-right) $p$-Yetter-Drinfeld module over $(H, \alpha)$ (or simply $\mathcal{YD}_{p}$-module)
 is an object $(M, \mu) \in \mathcal{\widetilde{H}}(\mathcal{M}_{k})$ such that
 \begin{enumerate}
   \item[(1)] $(M, \mu)$ is a left $(H_{p}, \alpha_{p})$-module (with notation $h_{p} \otimes m \mapsto h_{p} \cdot m$);
   \item[(2)] $(M, \mu, \rho^{M})$ is a right $(H, \alpha)$-module, i.e., there are a  family of linear maps
              $\rho^{M}=\{\rho^{M}_{p}\}_{p\in G}$, $\rho^{M}_{r}: M\longrightarrow M\otimes H_{r}, m\mapsto m_{(0)}\otimes m_{(1, r)}$
              in $\mathcal{\widetilde{H}}(\mathcal{M}_{k})$ such that
              \begin{eqnarray}
                \mu^{-1}(m_{(0)}) \otimes \Delta_{p, q}(m_{(1, pq)})
                   &=& m_{(0)(0)} \otimes (m_{(0)(1, p)} \otimes \alpha^{-1}_{q}(m_{(1, q)})), \label{1}\\
                \quad m_{(0)}\varepsilon(m_{(1, e)}) &=& \mu^{-1}(m),\label{2}
              \end{eqnarray}
              for all $m\in M$. Indeed, $\rho^{M} \in \mathcal{\widetilde{H}}(\mathcal{M}_{k})$ means for any $r \in G$
              \begin{eqnarray}
              \rho^{M}_{r}(\mu(m)) = \mu(m_{(0)}) \otimes \alpha_{r}(m_{(1, r)});
              \end{eqnarray}
   \item[(3)] the compatibility condition
              \begin{eqnarray}
              && h_{(1, p)} \cdot m_{(0)} \otimes h_{(2, r)} m_{(1, r)} \nonumber\\
              &=& \mu( (h_{(2, p)} \cdot \mu^{-1}(m))_{(0)} )\otimes
                     (h_{(2, p)} \cdot \mu^{-1}(m))_{(1, r)} \pi_{p^{-1}}(h_{(1, prp^{-1})}). \label{3}
              \end{eqnarray}
 \end{enumerate}

 Given two $\mathcal{YD}_{p}$-modules $(M, \mu, \rho^{M})$ and $(N, \nu, \rho^{N})$, a morphism of this two
 $\mathcal{YD}_{p}$-modules $f: (M, \mu, \rho^{M}) \longrightarrow (N, \nu, \rho^{N})$ is an $(H_{p}, \alpha_{p})$-linear
 map such that $f\circ\mu = \nu\circ f$ and for any $m\in M$ and $r \in G$
 \begin{eqnarray}
 \rho^{N}_{r} \circ f = (f \otimes Id_{H_{r}}) \circ \rho^{M}_{r}, \label{4}
 \end{eqnarray}
 i.e., $f(m)_{(0)} \otimes f(m)_{(1, r)} = f(m_{(0)}) \otimes m_{(1, r)}$.
 Then we have the category $\mathcal{YD}(H)_{p}$ of $\mathcal{YD}_{p}$-modules,
 the composition of morphisms of $\mathcal{YD}_{p}$-modules is the standard composition of the underlying linear maps.
 \\

 \textbf{Remark}
 (1) Equation (\ref{1}) exactly means
 \begin{eqnarray}
 \widetilde{a}_{M, H_{p}, H_{q}}(\rho^{M}_{p}\otimes Id_{H_{q}})\rho^{M}_{q} = (Id_{M} \otimes \Delta_{p, q}) \rho_{pq}.
 \end{eqnarray}
 Equation (\ref{3}) is equivalent to the following:
 \begin{eqnarray}
 && (h_{p}\cdot m)_{(0)} \otimes (h_{p}\cdot m)_{(1, r )} \nonumber \\
 &=& \alpha_{p}(h_{(12, p)}) \cdot m_{(0)}
 \otimes \alpha^{-1}_{r}(h_{(2, r)} m_{(1, r)}) S^{-1}_{r} \pi_{ p^{-1}}\alpha_{p r^{-1} p^{-1}}(h_{(11, p r^{-1} p^{-1})})  \label{5}
 \end{eqnarray}
 for all $h_{p} \in H_{p}$ and $m \in M \in\mathcal{YD}_{p}$.

 (2)  Let $(H=\bigoplus_{p\in G}H_{p}, \alpha = \{\alpha_{p}\}_{p\in G})$ be a monoidal Hom-Hopf $T$-coalgebra, then
 $(H, \alpha)$ is a left-right $p$-Yetter-Drinfeld module with the left $(H, \alpha)$-action
 \begin{eqnarray*}
 h\cdot g = \left\{
              \begin{array}{ll}
                (h_{(2, q)}\alpha^{-1}_{q}(g_{q}))S^{-1}\alpha_{q^{-1}}(h_{(1, q^{-1})}),
                   & \hbox{for} \quad h\in H_{e}, g\in H_{q}, \forall q\in G;\\
                0, & \hbox{otherwise.}
              \end{array}
            \right.
 \end{eqnarray*}
 and right $(H, \alpha)$-coaction by the Hom-comultiplication $\Delta = \{\Delta_{p, q}\}_{p, q\in G}$.
 $(k, Id_{k})$ is a left-right $p$-Yetter-Drinfeld module with the left $(H, \alpha)$-action
 \begin{eqnarray*}
 h\cdot \mathrm{k} = \left\{
              \begin{array}{ll}
               \varepsilon(h)\mathrm{k}, & \hbox{for} \quad h\in H_{e}, \mathrm{k}\in k;\\
                0, & \hbox{otherwise.}
              \end{array}
            \right.
 \end{eqnarray*}
 and right $(H, \alpha)$-coaction $\rho^{k}_{r}(\mathrm{k}) = \mathrm{k}\otimes 1_{H_{r}}$.
 \\

 \textbf{Proposition \thesection.2}
 If $(M, \mu, \rho^{M})\in \mathcal{YD}(H)_{p}$, $(N, \nu, \rho^{N})\in \mathcal{YD}(H)_{q}$,
 then $(M \otimes N, \mu \otimes \nu) \in \mathcal{YD}(H)_{pq}$ with the module and comodule structures as follows:
 \begin{eqnarray}
 h \cdot (m \otimes n)
 &=& h_{(1, p)} \cdot m \otimes h_{(2, q)} \cdot n,  \label{6} \\
 \rho^{M\otimes N}_{r}(m \otimes n)
 &=& m_{(0)} \otimes n_{(0)} \otimes n_{(1, r)}\pi_{q^{-1}}(m_{(1, qrq^{-1})}), \label{7}
 \end{eqnarray}
 where $m\in M, n\in N$, $h \in H_{pq}$.

 \emph{Proof} Obviously $V \otimes W$ is a left $(H_{pq}, \alpha_{pq})$-module, since for $h, k \in H_{pq}$,
 \begin{eqnarray*}
 \alpha_{pq}(h)\cdot (k \cdot (m \otimes n))
 &=& \alpha_{pq}(h)\cdot (k_{(1, p)} m \otimes k_{(2, q)} n) \\
 &=& (\alpha_{pq}(h))_{(1, p)} (k_{(1, p)} m) \otimes (\alpha_{pq}(h))_{(2, q)} (k_{(2, q)} n) \\
 &=& \alpha_{p}(h_{(1, p)}) (k_{(1, p)} m) \otimes \alpha_{q}(h_{(2, q)}) (k_{(2, q)} n) \\
 &=& (h_{(1, p)} k_{(1, p)}) \mu(m) \otimes (h_{(2, q)} k_{(2, q)})\nu(n) \\
 &=& (hk) \cdot ((\mu\otimes\nu)(m\otimes n)),\\
 1_{H_{pq}} \cdot (m \otimes n) &=& 1_{(1, p)} \cdot m \otimes 1_{(2, q)} \cdot n = (\mu\otimes\nu)(m\otimes n).
 \end{eqnarray*}
 In the following, we first check that coassociative condition holds.
 \begin{eqnarray*}
 && \widetilde{a}_{M\otimes N, H_{r_{1}}, H_{r_{2}}}(\rho^{M\otimes N}_{r_{1}}\otimes Id_{H_{r_{2}}})\rho^{M\otimes N}_{r_{2}} (m\otimes n) \\
 &=& \widetilde{a}_{M\otimes N, H_{r_{1}}, H_{r_{2}}}(\rho^{M\otimes N}_{r_{1}}\otimes Id_{H_{r_{2}}})
     (m_{(0)} \otimes n_{(0)} \otimes n_{(1, r_{2})}\pi_{q^{-1}}(m_{(1, qr_{2}q^{-1})}))  \\
 &=& \widetilde{a}_{M\otimes N, H_{r_{1}}, H_{r_{2}}} \\
  && \big((m_{(0)(0)} \otimes n_{(0)(0)} \otimes n_{(0)(1, r_{1})}\pi_{q^{-1}}(m_{(0)(1, qr_{1}q^{-1})}))
     \otimes n_{(1, r_{2})}\pi_{q^{-1}}(m_{(1, qr_{2}q^{-1})}) \big)  \\
 &=& (\mu\otimes\nu)(m_{(0)(0)} \otimes n_{(0)(0)}) \\
  && \otimes \big( n_{(0)(1, r_{1})}\pi_{q^{-1}}(m_{(0)(1, qr_{1}q^{-1})})
     \otimes \alpha^{-1}_{r_{2}}( n_{(1, r_{2})}\pi_{q^{-1}}(m_{(1, qr_{2}q^{-1})}) ) \big)  \\
 &=& \mu(m_{(0)(0)}) \otimes \nu(n_{(0)(0)}) \\
  && \otimes \big( n_{(0)(1, r_{1})}\pi_{q^{-1}}(m_{(0)(1, qr_{1}q^{-1})})
     \otimes \alpha^{-1}_{r_{2}}( n_{(1, r_{2})}) \pi_{q^{-1}}\alpha^{-1}_{r_{2}}(m_{(1, qr_{2}q^{-1})}) ) \big)  \\
 &\stackrel{(\ref{1})}{=}& (Id_{M\otimes N} \otimes \Delta_{r_{1}, r_{2}})
    \big(m_{(0)} \otimes n_{(0)} \otimes n_{(1, r_{1}r_{2})}\pi_{q^{-1}}(m_{(1, qr_{1}r_{2}q^{-1})}) \big) \\
 &=& (Id_{M\otimes N} \otimes \Delta_{r_{1}, r_{2}}) \rho^{M\otimes N}_{r_{1}r_{2}}(m\otimes n).
 \end{eqnarray*}

 Then the counitary condition (\ref{2}) is easy to get.

 Finally, we check the compatibility condition (\ref{3}) as follows.
 \begin{eqnarray*}
 && h_{(1, p q)} \cdot (m \otimes n)_{(0)} \otimes h_{(2, r)} (m \otimes n)_{(1, r)} \\
 &=& h_{(1, p q)} \cdot (m_{(0)} \otimes n_{(0)})
   \otimes h_{(2, r)} \big( n_{(1, r)}\pi_{q^{-1}}(m_{(1, q r q^{-1})}) \big)  \\
 &=& \big( h_{(11, p)} \cdot m_{(0)} \otimes h_{(12, q)} \cdot n_{(0)} \big)
   \otimes h_{(2, r)} \big( n_{(1, r)}\pi_{q^{-1}}(m_{(1, q r q^{-1})}) \big) \\
 &=& \big( h_{(11, p)} \cdot m_{(0)} \otimes h_{(12, q)} \cdot n_{(0)} \big)
   \otimes \big(\alpha^{-1}_{r}(h_{(2, r)}) n_{(1, r)}\big) \alpha_{r}\pi_{q^{-1}}(m_{(1, q r q^{-1})}) \\
 &=& \big( \alpha^{-1}_{p}(h_{(1, p)}) \cdot m_{(0)} \otimes \underline{h_{(21, q)} \cdot n_{(0)}} \big)
   \otimes \big(\underline{h_{(22, r)} n_{(1, r)}}\big) \alpha_{r}\pi_{q^{-1}}(m_{(1, q r q^{-1})}) \\
 &\stackrel{(\ref{1})}{=}&
   \big( \alpha^{-1}_{p}(h_{(1, p)}) \cdot m_{(0)} \otimes \nu[(h_{(22, q)} \cdot \nu^{-1}(n))_{(0)}] \big)\otimes \\
   &&\big( (h_{(22, q)} \cdot \nu^{-1}(n))_{(1, r)} \pi_{q^{-1}}(h_{(21, qrq^{-1})})\big) \alpha_{r}\pi_{q^{-1}}(m_{(1, q r q^{-1})})\\
 &=& \big( \alpha^{-1}_{p}(h_{(1, p)}) \cdot m_{(0)} \otimes \nu[(h_{(22, q)} \cdot \nu^{-1}(n))_{(0)}] \big)\otimes \\
   && \alpha[(h_{(22, q)} \cdot \nu^{-1}(n))_{(1, r)}]
   [\pi_{q^{-1}}(h_{(21, qrq^{-1})} m_{(1, q r q^{-1})})] \\
 &=& \big( \underline{h_{(11, p)} \cdot m_{(0)}} \otimes \nu[(\alpha^{-1}_{q}(h_{(2, q)}) \cdot \nu^{-1}(n))_{(0)}] \big)\otimes \\
   && \alpha[(\alpha^{-1}_{q}(h_{(2, q)}) \cdot \nu^{-1}(n))_{(1, r)}]
   \pi_{q^{-1}}[\underline{h_{(12, qrq^{-1})} m_{(1, q r q^{-1})}}] \\
 &\stackrel{(\ref{1})}{=}&
   \big( \mu[(h_{(12, p)} \cdot \mu^{-1}(m))_{(0)}] \otimes \nu[(\alpha^{-1}_{q}(h_{(2, q)}) \cdot \nu^{-1}(n))_{(0)}] \big)\otimes \\
   && \alpha[(\alpha^{-1}_{q}(h_{(2, q)}) \cdot \nu^{-1}(n))_{(1, r)}]
   \pi_{q^{-1}}[(h_{(12, p)} \cdot \mu^{-1}(m))_{(1, qrq^{-1})} \pi_{p^{-1}}(h_{(11, pqrq^{-1}p^{-1})})] \\
 &=& \big( \mu[(h_{(21, p)} \cdot \mu^{-1}(m))_{(0)}] \otimes \nu[(h_{(22, q)} \cdot \nu^{-1}(n))_{(0)}] \big)
   \otimes \alpha[(h_{(22, q)} \cdot \nu^{-1}(n))_{(1, r)}] \\
   && \pi_{q^{-1}}[(h_{(21, p)} \cdot \mu^{-1}(m))_{(1, qrq^{-1})} \pi_{p^{-1}}\alpha^{-1}_{pqrq^{-1}p^{-1}}(h_{(1, pqrq^{-1}p^{-1})})] \\
 &=& \big( \mu[(h_{(21, p)} \cdot \mu^{-1}(m))_{(0)}] \otimes \nu[(h_{(22, q)} \cdot \nu^{-1}(n))_{(0)}] \big)  \otimes \\
   && [(h_{(22, q)} \cdot \nu^{-1}(n))_{(1, r)}
   \pi_{q^{-1}}(h_{(21, p)} \cdot \mu^{-1}(m))_{(1, qrq^{-1})}] \pi_{q^{-1}p^{-1}}(h_{(1, pqrq^{-1}p^{-1})}) \\
 &=& (\mu\otimes\nu)[\underline{(h_{(21, p)} \cdot \mu^{-1}(m))_{(0)} \otimes (h_{(22, q)} \cdot \nu^{-1}(n))_{(0)}}]   \otimes \\
   && [\underline{(h_{(22, q)} \cdot \nu^{-1}(n))_{(1, r)}
   \pi_{q^{-1}}\big((h_{(21, p)} \cdot \mu^{-1}(m))_{(1, qrq^{-1})}\big)}] \pi_{q^{-1}p^{-1}}(h_{(1, pqrq^{-1}p^{-1})}) \\
 &\stackrel{(\ref{7})}{=}&
    (\mu\otimes\nu)[(h_{(21, p)} \cdot \mu^{-1}(m) \otimes h_{(22, q)} \cdot \nu^{-1}(n))_{(0)}]   \otimes \\
   && [(h_{(21, p)} \cdot \mu^{-1}(m) \otimes h_{(22, q)} \cdot \nu^{-1}(n))_{(1, r)}] \pi_{q^{-1}p^{-1}}(h_{(1, pqrq^{-1}p^{-1})}) \\
 &=& (\mu\otimes\nu)[(h_{(2, pq)} \cdot (\mu^{-1}(m) \otimes \nu^{-1}(n)))_{(0)}]   \otimes \\
   && [(h_{(2, pq)} \cdot (\mu^{-1}(m) \otimes \nu^{-1}(n)))_{(1, r)}] \pi_{q^{-1}p^{-1}}(h_{(1, pqrq^{-1}p^{-1})}).
 \end{eqnarray*}
 This completes the proof.
 $\hfill \Box$
 \\

 Following the left index notation shown in \cite{YW11} and \cite{Z04}, let $(M, \mu) \in \mathcal{YD}(H)_{p}$,
 the object ${}^{q}V$ has the same underling vector space as $V$.
 Given $v\in V$, we denote ${}^{q}v$ (${}^{q}f$) the corresponding element (morphism) in ${}^{q}V$.
 \\

 \textbf{Proposition \thesection.3}
 Let $(M, \mu, \rho^{M})\in \mathcal{YD}(H)_{p}$ and $q \in G$. Set ${}^{q}M=M$ as a vector space,
 with structures
 \begin{eqnarray}
   h_{q p q^{-1}} \rightharpoonup {}^{q}m &=& {}^{q}(\pi_{q^{-1}}(h_{q p q^{-1}}) \cdot m), \label{8} \\
   \rho^{{}^{q}M}_{r}({}^{q}m) &=& {}^{q}(m_{(0)}) \otimes \pi_{q}(m_{(1, q^{-1} r q)}) := ({}^{q}m)_{<0>} \otimes ({}^{q}m)_{<1,r>}, \label{9}
 \end{eqnarray}
 for any $m\in M$ and $a_{q p q^{-1}} \in H_{q p q^{-1}}$.
 Then ${}^{q}V \in \mathcal{YD}(H)_{q p q^{-1}}$.

 \emph{Proof}
 Obviously ${}^{q}V$ is a left $H_{q p q^{-1}}$-module,
 conditions (\ref{1}) and (\ref{2}) are also satisfied. In the following, we show that (\ref{5}) holds.
 \begin{eqnarray*}
 && (h_{q p q^{-1}} \rightharpoonup {}^{q}m)_{<0>} \otimes (h_{q p q^{-1}} \rightharpoonup {}^{q}m)_{<1, r>} \\
 &=& \rho^{{}^{q}M}_{r}(h_{q p q^{-1}} \rightharpoonup {}^{q}m)
  =  \rho^{{}^{q}M}_{r}({}^{q}(\pi_{q^{-1}}(h_{q p q^{-1}}) \cdot m)) \\
 &\stackrel{(\ref{9})}{=}&
    {}^{q}((\pi_{q^{-1}}(h_{q p q^{-1}}) \cdot m)_{(0)}) \otimes \pi_{q}[(\pi_{q^{-1}}(h_{q p q^{-1}}) \cdot m)_{(1, q^{-1}rq)}] \\
 &=& {}^{q}(\alpha_{p}\pi_{q^{-1}}(h_{(12, q p q^{-1})}) \cdot m_{(0)}) \otimes
    \pi_{q}[\alpha^{-1}_{q^{-1}rq}(\pi_{q^{-1}}(h_{(2, r)}) m_{(1, q^{-1}rq)}) \\
   &&  S^{-1}_{q^{-1}rq}\pi_{p^{-1}} \alpha_{ p q^{-1} r^{-1} q p^{-1}} \pi_{q^{-1}}(h_{(11, q p q^{-1} r^{-1} q p^{-1} q^{-1})})] \\
 &=& {}^{q}(\pi_{q^{-1}}(\alpha_{qpq^{-1}}(h_{(12, q p q^{-1})}) \cdot m_{(0)}) \otimes \\
   && \alpha^{-1}_{r}[h_{(2, r)} \pi_{q}(m_{(1, q^{-1}rq)}))]
    S^{-1}_{r}\pi_{qp^{-1}} \alpha_{ p q^{-1} r^{-1} q p^{-1}} \pi_{q^{-1}}(h_{(11, q p q^{-1} r^{-1} q p^{-1} q^{-1})})]\\
 &=& \alpha_{qpq^{-1}}(h_{(12, q p q^{-1})}) \rightharpoonup {}^{q}(m_{(0)}) \otimes \\
   && \alpha^{-1}_{r}[h_{(2, q p q^{-1})} \pi_{q}(m_{(1, q^{-1}rq)}))]
    S^{-1}_{r}\pi_{qp^{-1}q^{-1}} \alpha_{q p q^{-1} r q p^{-1} q^{-1}}(h_{(11, q p q^{-1} r q p^{-1} q^{-1})}) \\
 &=& \alpha_{qpq^{-1}}(h_{(12, q p q^{-1})}) \rightharpoonup ({}^{q}m)_{<0>} \otimes \\
   && \alpha^{-1}_{r}[h_{(2, q p q^{-1})} ({}^{q}m)_{<1, r>} )]
    S^{-1}_{r}\pi_{qp^{-1}q^{-1}} \alpha_{q p q^{-1} r q p^{-1} q^{-1}}(h_{(11, q p q^{-1} r q p^{-1} q^{-1})})
 \end{eqnarray*}
 This completes the proof.
 $\hfill \Box$
 \\

 \textbf{Remark}
 Let $(M, \mu, \rho^{M})\in \mathcal{YD}(H)_{p}$ and $(N, \nu, \rho^{N})\in \mathcal{YD}(H)_{q}$,
 then we have ${}^{s t}M={}^{s}({}^{t}M)$
 as an object in $\mathcal{YD}(H)_{s t p t^{-1} s^{-1}}$, and
 ${}^{s}(M \otimes N) = {}^{s}M \otimes {}^{s}N$
 as an object in $\mathcal{YD}(H)_{s p q s^{-1}}$.
 \\

 \textbf{Proposition \thesection.4}
 Let $(M, \mu, \rho^{M})\in \mathcal{YD}(H)_{p}$ and $(N, \nu, \rho^{N})\in \mathcal{YD}(H)_{q}$.
 Set ${}^{M}N={}^{p}N$ as an object in $\mathcal{YD}(H)_{p q p^{-1}}$.
 Define the map
 \begin{eqnarray}
 c_{M, N} &:& M\otimes N \longrightarrow {}^{M}N \otimes M, \nonumber \\
 c_{M, N}(m \otimes n) &=& {}^{p}(S_{q^{-1}}(m_{(1, q^{-1})}) \cdot \nu^{-1}(n)) \otimes \mu(m_{(0)}). \label{10}
 \end{eqnarray}
 Then $c_{M, N}$ is $(H, \alpha)$-linear, $(H, \alpha)$-colinear and satisfies the the Hexagon axiom in \cite{K95}
 \begin{eqnarray}
  (Id_{{}^{M}N} \otimes c_{M, X})\circ \widetilde{a}_{{}^{M}N, M, X} \circ (c_{M, N}\otimes Id_{X})
    = \widetilde{a}_{{}^{M}N, {}^{M}X, M} \circ c_{M, N\otimes X}\circ \widetilde{a}_{M, N, X} \label{11} \\
  \widetilde{a}^{-1}_{{}^{M\otimes N}X, M, N} \circ c_{M\otimes N, X}\circ \widetilde{a}^{-1}_{M, N, X}
 = (c_{M, {}^{N}X} \otimes Id_{N}) \circ \widetilde{a}^{-1}_{M, {}^{N}X, N}\circ (Id_{M} \otimes c_{N, X}) \label{12}
 \end{eqnarray}
 for $(X, \eta, \rho^{X}) \in \mathcal{YD}(H)_{s}$. Moreover, $c_{{}^{s}M, {}^{s}N} = {}^{s}(\cdot) \circ c_{M, N}$.

 \emph{Proof}
 Firstly, we prove that $c_{M, N}$ is $(H, \alpha)$-linear. We compute
 \begin{eqnarray*}
 && c_{M, N}(h_{p q} \cdot (m \otimes n)) \\
 &=& c_{M, N}(h_{(1, p)} \cdot m \otimes h_{(2, q)} \cdot n) \\
 &=& {}^{p}(S_{q^{-1}}((h_{(1, p)} \cdot m)_{(1, q^{-1})}) \cdot \nu^{-1}(h_{(2, q)} \cdot n))
   \otimes \mu((h_{(1, p)} \cdot m)_{(0)}) \\
 &=& {}^{p}\big( S_{q^{-1}}[\alpha^{-1}_{q^{-1}}(h_{(12, q^{-1})} m_{(1, q^{-1})})
     S^{-1}_{q^{-1}}\pi_{p^{-1}}\alpha_{pqp^{-1}}(h_{(111, pqp^{-1})})] \cdot \nu^{-1}(h_{(2, q)} \cdot n) \big) \\
  && \otimes  \mu( \alpha_{p}(h_{(112, p)}) \cdot m_{(0)}) \\
 &=& {}^{p}\big([\pi_{p^{-1}}\alpha_{pqp^{-1}}(h_{(111, pqp^{-1})})
     S_{q^{-1}}\alpha^{-1}_{q^{-1}}(h_{(12, q^{-1})} m_{(1, q^{-1})})] \cdot( \alpha^{-1}_{q}(h_{(2, q)}) \cdot \nu^{-1}(n) ) \big) \\
  && \otimes  \mu(\alpha_{p}(h_{(112, p)}) \cdot m_{(0)}) \\
 &=& {}^{p}\Big([ \big(\pi_{p^{-1}}(h_{(111, pqp^{-1})})
     S_{q^{-1}}\alpha^{-2}_{q^{-1}}(h_{(12, q^{-1})} m_{(1, q^{-1})}) \big) \alpha^{-1}_{q}(h_{(2, q)})] \cdot n ) \Big) \\
  && \otimes  \mu(\alpha_{p}(h_{(112, p)}) \cdot m_{(0)}) \\
 &=& {}^{p}\Big([\alpha_{q} \pi_{p^{-1}}(h_{(111, pqp^{-1})})
     \big(S_{q^{-1}}\alpha^{-2}_{q^{-1}}(h_{(12, q^{-1})} m_{(1, q^{-1})}) \alpha^{-2}_{q}(h_{(2, q)}) \big) ] \cdot n ) \Big) \\
  && \otimes  \mu( \alpha_{p}(h_{(112, p)}) \cdot m_{(0)}) \\
 &=& {}^{p}\Big([\alpha_{q} \pi_{p^{-1}}(h_{(111, pqp^{-1})})
     \big( [S_{q^{-1}}\alpha^{-2}_{q^{-1}}(m_{(1, q^{-1})}) S_{q^{-1}}\alpha^{-2}_{q^{-1}}(h_{(12, q^{-1})})]
     \alpha^{-2}_{q}(h_{(2, q)}) \big) ] \cdot n ) \Big) \\
  && \otimes  \mu( \alpha_{p}(h_{(112, p)}) \cdot m_{(0)}) \\
 &=& {}^{p}\Big([\alpha_{q} \pi_{p^{-1}}(h_{(111, pqp^{-1})})
     \big( S_{q^{-1}}\alpha^{-1}_{q^{-1}}(m_{(1, q^{-1})}) [S_{q^{-1}}\alpha^{-2}_{q^{-1}}(h_{(12, q^{-1})})
     \alpha^{-3}_{q}(h_{(2, q)})\big) ] \cdot n ) \Big) \\
  && \otimes \mu( \alpha_{p}(h_{(112, p)}) \cdot m_{(0)}) \\
 &=& {}^{p}\Big([ \pi_{p^{-1}}\alpha^{-1}_{pqp^{-1}}(h_{(1, pqp^{-1})})
     \big( S_{q^{-1}}\alpha^{-1}_{q^{-1}}(m_{(1, q^{-1})}) [S_{q^{-1}}\alpha^{-1}_{q^{-1}}(h_{(221, q^{-1})})
     \alpha^{-1}_{q}(h_{(222, q)})] \big) ] \cdot n ) \Big) \\
  && \otimes \mu( h_{(21, p)} \cdot m_{(0)}) \\
 &=& {}^{p}\Big([ \pi_{p^{-1}}\alpha^{-1}_{pqp^{-1}}(h_{(1, pqp^{-1})})
     \big( S_{q^{-1}}\alpha^{-1}_{q^{-1}}(m_{(1, q^{-1})}) \varepsilon(h_{(22, e)} 1_{H_{q}})\big) ] \cdot n ) \Big) \\
  && \otimes \mu( h_{(21, p)} \cdot m_{(0)}) \\
 &=& {}^{p}\Big([ \pi_{p^{-1}}\alpha^{-1}_{pqp^{-1}}(h_{(1, pqp^{-1})}) S_{q^{-1}}(m_{(1, q^{-1})}) ] \cdot n ) \Big)
   \otimes \mu(\alpha^{-1}_{p}(h_{(2, p)}) \cdot m_{(0)}),
 \end{eqnarray*}
 and
 \begin{eqnarray*}
 && h_{p q} \cdot c_{M, N}(m \otimes n) \\
 &=& h_{p q} \cdot ( {}^{p}(S_{q^{-1}}(m_{(1, q^{-1})}) \cdot \nu^{-1}(n)) \otimes \mu(m_{(0)}) ) \\
 &=& h_{(1, p q p^{-1})} \cdot {}^{p}(S_{q^{-1}}(m_{(1, q^{-1})}) \cdot \nu^{-1}(n))
   \otimes h_{(2, p)} \cdot \mu(m_{(0)}) \\
 &\stackrel{(\ref{8})}{=}&
   {}^{p}(\pi_{p^{-1}}(h_{(1, p q p^{-1})}) \cdot (S_{q^{-1}}(m_{(1, q^{-1})}) \cdot \nu^{-1}(n)) )
   \otimes h_{(2, p)} \cdot \mu(m_{(0)})\\
 &=& {}^{p}\big( (\alpha_{q} \pi_{p^{-1}}(h_{(1, p q p^{-1})}) S_{q^{-1}}(m_{(1, q^{-1})})) \cdot n \big)
   \otimes h_{(2, p)} \cdot \mu(m_{(0)}).
 \end{eqnarray*}

 Secondly, we can check that $c_{M, N}$ is $(H, \alpha)$-colinear, i.e., $\rho^{{}^{M}N \otimes M}_{r} \circ c_{M, N}
 = (c_{M, N} \otimes Id_{H_{r}}) \circ \rho^{M \otimes N}_{r}$.
 In fact,
 \begin{eqnarray*}
 && \rho^{{}^{M}N \otimes M}_{r} \circ c_{M, N} (m \otimes n)  \\
 &=& \rho^{{}^{M}N \otimes M}_{r} \Big({}^{p}(S_{q^{-1}}(m_{(1, q^{-1})}) \cdot \nu^{-1}(n)) \otimes \mu(m_{(0)}) \Big) \\
 &=& [{}^{p}(S_{q^{-1}}(m_{(1, q^{-1})}) \cdot \nu^{-1}(n))]_{<0>} \otimes \mu(m_{(0)})_{(0)}  \otimes \\
  && \mu(m_{(0)})_{(1, r)} \pi_{p^{-1}}\Big( [{}^{p}(S_{q^{-1}}(m_{(1, q^{-1})}) \cdot \nu^{-1}(n))]_{<1, prp^{-1}>} \Big) \\
 &=& {}^{p}(\underline{\big(S_{q^{-1}}(m_{(1, q^{-1})}) \cdot \nu^{-1}(n)\big)_{(0)}} ) \otimes \mu(m_{(0)})_{(0)}  \otimes \\
  && \mu(m_{(0)})_{(1, r)} \pi_{p^{-1}}\Big(\pi_{p} [\underline{\big(S_{q^{-1}}(m_{(1, q^{-1})}) \cdot \nu^{-1}(n) \big)_{(1, r)}}] \Big) \\
 &\stackrel{(\ref{5})}{=}&
     {}^{p}(\alpha_{q} S_{q^{-1}}(m_{(121, q^{-1})}) \cdot \nu^{-1}(n_{(0)}) ) \otimes \mu(m_{(0)(0)})  \otimes \\
  && \alpha(m_{(0)(1, r)})
    \Big( \alpha^{-1}_{r}[S_{r^{-1}}(m_{(11,r)}) \alpha^{-1}_{r}(n_{(1, r)})] \pi_{q^{-1}}\alpha_{qrq^{-1}} (m_{(122, qrq^{-1})}) \Big) \\
 &=& {}^{p}(\alpha_{q} S_{q^{-1}}(m_{(121, q^{-1})}) \cdot \nu^{-1}(n_{(0)}) ) \otimes \mu(m_{(0)(0)})  \otimes \\
  && \Big(m_{(0)(1, r)} (\alpha^{-1}_{r}[S_{r^{-1}}(m_{(11,r^{-1})}) \alpha^{-1}_{r}(n_{(1, r)})]) \Big)
     \pi_{q^{-1}}\alpha^{2}_{qrq^{-1}} (m_{(122, qrq^{-1})}) \\
 &=& {}^{p}(\alpha_{q} S_{q^{-1}}(m_{(121, q^{-1})}) \cdot \nu^{-1}(n_{(0)}) ) \otimes \mu(m_{(0)(0)})  \otimes \\
  && \Big( [\alpha^{-1}_{r}(m_{(0)(1, r)}) \alpha^{-1}_{r}S_{r^{-1}}(m_{(11,r^{-1})})] \alpha^{-1}_{r}(n_{(1, r)}) \Big)
     \pi_{q^{-1}}\alpha^{2}_{qrq^{-1}} (m_{(122, qrq^{-1})}) \\
 &=& {}^{p}(\alpha_{q} S_{q^{-1}}(m_{(121, q^{-1})}) \cdot \nu^{-1}(n_{(0)}) ) \otimes m_{(0)}  \otimes \\
  && \Big( [\alpha^{-1}_{r}\alpha_{r}(m_{(111, r)}) \alpha^{-1}_{r}S_{r^{-1}}\alpha_{r^{-1}}(m_{(112,r^{-1})})]
     \alpha^{-1}_{r}(n_{(1, r)}) \Big)
     \pi_{q^{-1}}\alpha^{2}_{qrq^{-1}} (m_{(122, qrq^{-1})}) \\
 &=& {}^{p}(\alpha_{q} S_{q^{-1}}(m_{(121, q^{-1})}) \cdot \nu^{-1}(n_{(0)}) ) \otimes m_{(0)}  \otimes
     n_{(1, r)} \varepsilon(m_{11, e})  \pi_{q^{-1}}\alpha^{2}_{qrq^{-1}} (m_{(122, qrq^{-1})}) \\
 &=& {}^{p}(\alpha_{q} S_{q^{-1}}\alpha_{q^{-1}}(m_{(0)(1, q^{-1})} \cdot \nu^{-1}(n_{(0)}) ) \otimes \mu(m_{(0)(0)})  \otimes\\
  && n_{(1, r)} \pi_{q^{-1}}\alpha^{2}_{qrq^{-1}} \alpha^{-2}_{qrq^{-1}} (m_{(1, qrq^{-1})}) \\
 &=& {}^{p}( S_{q^{-1}}(m_{(0)(1, q^{-1})} \cdot \nu^{-1}(n_{(0)}) ) \otimes \mu(m_{(0)(0)})  \otimes
    n_{(1, r)} \pi_{q^{-1}} (m_{(1, qrq^{-1})}),
 \end{eqnarray*}
 and
 \begin{eqnarray*}
 && (c_{M, N} \otimes Id_{H_{r}}) \circ \rho^{M \otimes N}_{r}(m\otimes n) \\
 &=& (c_{M, N} \otimes Id_{H_{r}}) (m_{(0)} \otimes n_{(0)} \otimes n_{(1, r)}\pi_{q^{-1}}(m_{(1, qrq^{-1})})) \\
 &=& ({}^{p}(S_{q^{-1}}(m_{(0)(1, q^{-1})}) \cdot \nu^{-1}(n_{(0)})) \otimes \mu(m_{(0)(0)}))
    \otimes n_{(1, r)}\pi_{q^{-1}}(m_{(1, qrq^{-1})})).
 \end{eqnarray*}

 Thirdly, $c_{V, W}$ satisfies the conditions (\ref{11}) and (\ref{12}). Here we only check (\ref{11}), condition (\ref{12}) is similar.
 \begin{eqnarray*}
 && (Id_{{}^{M}N} \otimes c_{M, X})\circ \widetilde{a}_{{}^{M}N, M, X} \circ (c_{M, N}\otimes Id_{X})((m\otimes n) \otimes x) \\
 &=& (Id_{{}^{M}N} \otimes c_{M, X})\circ \widetilde{a}_{{}^{M}N, M, X}
     \Big([{}^{p}(S_{q^{-1}}(m_{(1, q^{-1})}) \cdot \nu^{-1}(n)) \otimes \mu(m_{(0)})] \otimes x \Big)\\
 &=& (Id_{{}^{M}N} \otimes c_{M, X})
     \Big(\nu \big({}^{p}(S_{q^{-1}}(m_{(1, q^{-1})}) \cdot \nu^{-1}(n)) \big) \otimes [\mu(m_{(0)}) \otimes \eta^{-1}(x)] \Big)\\
 &=& \nu \big({}^{p}(S_{q^{-1}}(m_{(1, q^{-1})}) \cdot \nu^{-1}(n)) \big) \otimes
     \Big( {}^{p}(S_{s^{-1}}(\mu(m_{(0)})_{(1, s^{-1})}) \cdot \eta^{-2}(x)) \otimes \mu(\mu(m_{(0)})_{(0)})  \Big) \\
 &=& \nu \big({}^{p}(S_{q^{-1}}(m_{(1, q^{-1})}) \cdot \nu^{-1}(n)) \big) \otimes
     \Big( {}^{p}(S_{s^{-1}}(\alpha_{s^{-1}}(m_{(0)(1, s^{-1})})) \cdot \eta^{-2}(x)) \otimes \mu^{2}(m_{(0)(0)})  \Big) \\
 &=& \nu \big({}^{p}(S_{q^{-1}}\alpha_{q^{-1}}(m_{(12, q^{-1})}) \cdot \nu^{-1}(n)) \big) \otimes
     \Big( {}^{p}(S_{s^{-1}}(\alpha_{s^{-1}}(m_{(11, s^{-1})})) \cdot \eta^{-2}(x)) \otimes \mu(m_{(0)})  \Big) \\
 &=&{}^{p}\big( S_{q^{-1}}\alpha^{2}_{q^{-1}}(m_{(12, q^{-1})}) \cdot n \big) \otimes
     \Big( {}^{p}(S_{s^{-1}}(\alpha_{s^{-1}}(m_{(11, s^{-1})})) \cdot \eta^{-2}(x)) \otimes \mu(m_{(0)})  \Big),
 \end{eqnarray*}
 and
 \begin{eqnarray*}
 && \widetilde{a}_{{}^{M}N, {}^{M}X, M} \circ c_{M, N\otimes X}\circ \widetilde{a}_{M, N, X} ((m\otimes n) \otimes x) \\
 &=& \widetilde{a}_{{}^{M}N, {}^{M}X, M} \circ c_{M, N\otimes X} (\mu(m) \otimes (n \otimes \eta^{-1}(x))) \\
 &=& \widetilde{a}_{{}^{M}N, {}^{M}X, M}
     \Big[ {}^{p}\Big( S_{(qs)^{-1}}(\mu(m)_{(1, (qs)^{-1})}) \cdot (\nu^{-1}(n) \otimes \eta^{-2}(x)) \Big) \otimes \mu(\mu(m)_{(0)}) \Big] \\
 &=& \widetilde{a}_{{}^{M}N, {}^{M}X, M}
     \Big[ {}^{p}\Big( S_{(qs)^{-1}}\alpha_{(qs)^{-1}} (m_{(1, (qs)^{-1})}) \cdot (\nu^{-1}(n) \otimes \eta^{-2}(x)) \Big)
     \otimes \mu^{2}(m_{(0)}) \Big] \\
 &=& \widetilde{a}_{{}^{M}N, {}^{M}X, M} \\
   && \Big[ {}^{p}\Big( S_{q^{-1}}\alpha_{q^{-1}} (m_{(12, q^{-1})}) \cdot \nu^{-1}(n) \otimes
     S_{s^{-1}}\alpha_{s^{-1}} (m_{(11, s^{-1})}) \cdot\eta^{-2}(x) \Big)
     \otimes \mu^{2}(m_{(0)}) \Big] \\
 &=& \nu( {}^{p}\Big( S_{q^{-1}}\alpha_{q^{-1}} (m_{(12, q^{-1})}) \cdot \nu^{-1}(n) \Big) ) \otimes \Big[
     {}^{p}\Big( S_{s^{-1}}\alpha_{s^{-1}} (m_{(11, s^{-1})}) \cdot\eta^{-2}(x) \Big)
     \otimes \mu(m_{(0)}) \Big] \\
 &=& {}^{p}\Big( S_{q^{-1}}\alpha^{2}_{q^{-1}} (m_{(12, q^{-1})}) \cdot n \Big)  \otimes \Big[
     {}^{p}\Big( S_{s^{-1}}\alpha_{s^{-1}} (m_{(11, s^{-1})}) \cdot\eta^{-2}(x) \Big)
     \otimes \mu(m_{(0)}) \Big].
 \end{eqnarray*}

 Finally, we check the last condition $c_{{}^{s}M, {}^{s}N} = {}^{s}(\cdot) \circ c_{M, N}$. Indeed,
 \begin{eqnarray*}
 && c_{{}^{s}M, {}^{s}N} ({}^{s}m \otimes {}^{s}n) \\
 &=& {}^{sps^{-1}}(S_{sq^{-1}s^{-1}}(({}^{s}m)_{<1, sq^{-1}s^{-1}>}) \rightharpoonup \nu^{-1}({}^{s}n)) \otimes \mu(({}^{s}m)_{<0>}) \\
 &=& {}^{sps^{-1}}(S_{sq^{-1}s^{-1}}(\pi_{s}(m_{(1, q^{-1})})) \rightharpoonup {}^{s}(\nu^{-1}(n))) \otimes \mu({}^{s}(m_{(0)})) \\
 &=& {}^{sps^{-1}}\Big({}^{s}[\pi_{s^{-1}}S_{sq^{-1}s^{-1}}(\pi_{s}(m_{(1, q^{-1})})) \cdot \nu^{-1}(n)] \Big)
     \otimes \mu({}^{s}(m_{(0)})) \\
 &=& {}^{sp}\Big(S_{q^{-1}}(m_{(1, q^{-1})}) \cdot \nu^{-1}(n) \Big) \otimes {}^{s}(\mu(m_{(0)})) \\
 &=& {}^{s}\Big( {}^{p}\big(S_{q^{-1}}(m_{(1, q^{-1})}) \cdot \nu^{-1}(n) \big) \otimes \mu(m_{(0)}) \Big) \\
 &=& {}^{s}(\cdot) \circ c_{M, N}(m \otimes n).
 \end{eqnarray*}
 This completes the proof.
 $\hfill \Box$
 \\

 Define $\mathcal{YD}(H)$ as the disjoint union of all $\mathcal{YD}(H)_{p}$ with $ p \in G$.
 If we endow $\mathcal{YD}(H)$ with tensor product as in Proposition \thesection.2, then it becomes a monoidal
 category. The unit is $(k, Id_{k})$ with trivial structures.

 The group homomorphism $\phi: G\longrightarrow aut(\mathcal{YD}(H)), p \mapsto \phi_{p}$ is given on
 components as
 \begin{eqnarray*}
 \phi_{p}: \mathcal{YD}(H)_{q} \longrightarrow \mathcal{YD}(H)_{p q p^{-1}},
 \phi_{p}(N)={}^{p}N,
 \end{eqnarray*}
 where the functor $\phi_{p}$ act as follows:
 given a morphism $f: (M, \mu, \rho^{M})\longrightarrow (N, \nu, \rho^{N})$, for any $m\in M$, we set
 $({}^{p}f)({}^{p}m)={}^{p}(f(m))$.

 The braiding in $\mathcal{YD}(H)$ is given by the family $\{c_{M, N} \}$.
 As a consequence of the above results, we obtain the main result of this paper:
 \\

 \textbf{Theorem \thesection.5}
 For a monoidal Hom-Hopf $T$-coalgebra $(H, \alpha)$, $\mathcal{YD}(H)$ is a braided $T$-category over group $G$.
 \\

 The braiding in $\mathcal{YD}(H)$ is invertible, the inverse is given as follows.

 \textbf{Proposition \thesection.6}
 Let $(M, \mu, \rho^{M})\in \mathcal{YD}(H)_{p}$ and $(N, \nu, \rho^{N})\in \mathcal{YD}(H)_{q}$.
 Then
 \begin{eqnarray}
 c_{M, N}^{-1}: {}^{M}N \otimes M \longrightarrow M\otimes N, \qquad
 c_{M, N}^{-1}({}^{p}n \otimes m) = \mu(m_{(0)}) \otimes m_{(1, q)} \cdot \nu^{-1}(n). \label{13}
 \end{eqnarray}

 \emph{Proof} It is straightforward.
 $\hfill \Box$

\section*{Acknowledgements}

 The authors would like to thank the referee for his/her valuable comments.
 The work was partially supported by the National Natural Science Foundation of China (Grant No. 11226070),
 the Fundamental Research Fund for the Central Universities(Grant No.KYZ201322)
 and the NJAUF (No. LXY201201019, LXYQ201201103).

\end {document}